\documentclass[12pt]{article} 

\usepackage{color}
\usepackage{latexsym}       
\usepackage{amssymb}
\usepackage{amscd}   
\def \diag{\mbox{\rm diag}} 
\usepackage{amssymb,amsthm,amsmath}
\usepackage[T1]{fontenc}
\newtheorem{theorem}{Theorem}[section]
\newcommand{\theo}{\begin{theorem}}
\newcommand{\enth}{\end{theorem}}

\newtheorem{proposition}[theorem]{Proposition}
\newcommand{\prop}{\begin{proposition}}
\newcommand{\enpro}{\end{proposition}}
\newcommand{\nothing}[1]{}

\newtheorem{Remark}[theorem]{Remark}
\newcommand{\remark}{\begin{Remark}}
\newcommand{\enremark}{\end{Remark}}

\newtheorem{Lemma}[theorem]{Lemma}

\newcommand{\ydgen}[5]
{
  \ifnum#2=3
\ifnum#1=1
    \left[y_{17}\right]_{\zeta^{#1}}\!d_{#3,#4,#5}
\fi
\ifnum#1=2
    \left[y_{16}\right]_{\zeta^{#1}}\!d_{#3,#4,#5}
\fi
\ifnum#1=4
    \left[y_{15}\right]_{\zeta^{#1}}\!d_{#3,#4,#5}
\fi

  \else
    \left[y_{#2}\right]_{\zeta^{#1}}\!d_{#3,#4,#5}
  \fi
}
\newcommand{\ydDj}[4]
{
  \ifnum#2=3
\ifnum#1=1
    \left[y_{17}\right]_{\zeta^{#1}}\!\left[D_{#3}\right]_{\zeta^{#4}}
\fi
\ifnum#1=2
    \left[y_{16}\right]_{\zeta^{#1}}\!\left[D_{#3}\right]_{\zeta^{#4}}
\fi
\ifnum#1=4
    \left[y_{15}\right]_{\zeta^{#1}}\!\left[D_{#3}\right]_{\zeta^{#4}}
\fi

  \else
    \left[y_{#2}\right]_{\zeta^{#1}}\!\left[D_{#3}\right]_{\zeta^{#4}}
  \fi
}
\newcommand{\ydgenfrstp}[5]
{
  \ifnum#2=2
\ifnum#1=13
    \left[y_{2}\right]_{\nu^{#1}}\!d_{#3,#4,#5}
\fi  
\else
    \left[y_{#2}\right]_{\nu^{#1}}\!d_{#3,#4,#5}
  \fi
}

\title{A canonical curve of genus $17$}
\author{Israel Moreno-Mej\'{\i}a}
\begin{document}
\maketitle
\begin{abstract}
We compute equations for a Hurwitz curve of genus 17 and we conclude that the canonical ideal of any Hurwitz curve of genus 14 or 17 is generated by quadrics.\footnote{2000 MSC 14H37,14H45, (30F10,20H10)}
 \end{abstract}
\section{Introduction}

The curve alluded in the title is one of the 2 Hurwitz curves of genus 17 and the aim of the paper is   to compute equations  for its  canonical embedding  in $\mathbb{P}^{16}$. 
Great deal of the work is about finding the quadrics $I_2$ in the curve's homogeneous ideal $I$ (see Theorem \ref{tbasI2}) and the rest is to verify that these quadrics generate $I$.
The results of this paper were originally meant to be part of \cite{hcg14} however some of the computations involved proved to be very tough and we could not reach any conclusion until long afterwards.
So, to compute the  space of quadrics in $I$ we follow the method used in \cite{hcg14}  for the case of the Hurwitz curves of genus 14; that is, we study the action of the curve's automorphism group $G$ on the space $H^{0}(X_{},K_{X_{}}^{\otimes d})$ of global sections of tensor powers of the canonical sheaf (see \S \ref{sectionquadrics}). Then using  a matrix representation of $G$ for this action one tries to compute distinguished subspaces of  $I_2$ whose sum is $I_2$, see \S\ref{secpartW4}, \S\ref{secpartW10} and \S\ref{secpartW11}.
We did not know the matrix representation of $G$ for the action on the space $H^{0}(X_{},K_{X_{}}^{})$  until Vahid Dabbaghian implemented in his package Repsn  of GAP \cite{GAP4}  an algorithm  to compute complex matrix representations of finite groups that worked for $G$ (see \cite{repsn}). Then the computations with some of the polynomials that we obtained turned out to be very slow and we were unable to work out the spaces corresponding to   \S\ref{secpartW10} and \S\ref{secpartW11} until we used ParGap \cite{Coo99} to performe some  parallel calculations to find the conjugating matrices  (\ref{cnjmtrx}) and (\ref{cnjmtrx2}). We did not need such matrices in \cite{hcg14} in the genus 14 case because the eigenspaces involved were 1-dimensional.

That $I$ is generated by quadrics would follow at once from  \emph{ Petri's Analysis of the Canonical Ideal} (see for instance \cite{Petrioriginal},\cite{petri1},\cite{mumcj} or \cite{petri2}) if one knew that the curve is not trigonal and although  this is the case for  a general curve of genus $g>4$ a proof is required for this specific curve (we also check the genus $14$ case). We used the quadrics to test the non-trigonality of the curve, see Theorem \ref{mainthm} for the genus 17 Hurwitz curve and Theorem \ref{mainthm2} for the genus 14  case. Having the quadrics in the ideal available to do calculations this seemed to be a quick way to know the answer,  however, this also turned out to be a very slow and long computation (at least with the algorithms that we used and I did not try to parallelize it, the tables in \S\ref{sectcubicforms}, \S\ref{sectgenus14} are  meant to help speed up verification) and a simpler method would be desirable. Perhaps using a different interpretation of the curve one could apply tools similar to those  used to determine the gonanality of modular curves (see for instance \cite{TMC}, \cite{TQoMC} ). Nevertheless using the quadrics in the ideal to check the trigonality of modular curves is also among 
those tools (see \S 2 in \cite{TMC} or \S 2.2 in \cite{DEoMC}).

One can convince oneself from the calculation of the quadrics that, it is possible to obtain equations for the remaining genus 17 and genus 14 Hurwitz curves by applying appropriate elements of the absolute Galois group  $Gal(\bar{\mathbb{Q}},\mathbb{Q})$ to the coefficients of the equations of the curves considered here. In the genus 14 case, this had already been observed in \cite{streit}. It follows that the canonical ideal of any  genus 14 or genus 17 Hurwitz curve is generated by quadrics.

Very little is known about the genus 17 Hurwitz curves and the information available to us about their automorphism group comes from \cite{conder} and \cite{sinkov}. One could highlight the fact that there is only one Hurwitz group $G$ of genus $17$ and that it can be realized as the automorphism group of 2 non isomorphic genus 17 curves $X_{1}$ and $X_{2}$  (see Section \ref{sctnx1} for the definition of $X_{i}$). Unlike the Hurwitz groups of genus $3$,$7$ and $14$ our group $G$ is not simple and has a normal subgroup $N$   isomorphic to $(\mathbb{Z}/2\mathbb{Z})^{3}$. One has that $G/N$ is $PSL(2,\mathbb{F}_{7})$ and $X_{1}/N\cong X_{2}/N$ is the Hurwitz curve of genus $3$. This of course plays a role in our calculation, although in a very subtle way, it helped us when we tried to find a point of  the curve in $\S\ref{secpartW4}$ which was crucial to find some of the equations.  

We remark that it is not known a way to represent   a general curve of genus $g$ by means of polynomial equations, except for small values of $g$ (see \cite{moduliofcurves} Chap.6 ); and  that although for some genera  it is possible to give equations of special curves,   it appears that for the genus 14 and genus 17  cases the  Hurwitz curves  are the only ones so far for which equations defining them can be computed.

The problem of computing equations for the Hurwitz curves of genus 14 and 17 had been an open problem for quite a while, see for instance  the Problem 7 in \cite{shaska} and  $\S$8 in \cite{macbeathproblem}. Unfortunately most of the  polynomials  that we found have coefficients that are impossible to display here because of their size, however, producing these equations with GAP o ParGap is rather quick even in the genus 17 case (using the formulae in $\S$\ref{basI2})  if one has at hand the conjugating matrices (\ref{cnjmtrx}) and (\ref{cnjmtrx2}) which can be obtained from me. For doing calculations with them it would be better to use a 64--bit version of GAP or ParGap and to have available upto 8 gigabytes of RAM memory.

\section{Notation}\label{sectionnotation}
We keep our notation as in \cite{hcg14}. Given a curve $X$ and $h\in$Aut($X$) 
with fixed point set $fix(h)=\{p_1,\ldots,p_u\}$ then 
\[d_{X}h=(dh_{p_{1}},\ldots, dh_{p_{u}}),\]
where $dh_{p_{i}}$ is the  eigenvalue of $h$ acting on 
$K_{X,p_{i}}$, where $K_{X,p_{i}}$ is the fibre of the canonical line bundle of $X$ at $p_{i}$.\\
If $G$ is a finite group  then $\textsf{W}_{1},\ldots, \textsf{W}_{s}$ denote its complex irreducible representations 
and  $\chi_{1},\ldots \chi_{s}$ denote the corresponding characters. Given   a $\mathbb{C}G$-module $M$  we write $M_{\textsf{W}_{i}}$ for the sum of all $\mathbb{C}G$-submodules of $M$ isomorphic to $W_{i}$. The projection
$M\rightarrow M_{\textsf{W}_{i}}$ is denoted by $\pi_{\textsf{W}_{i}}$. 
For $h\in G$, $\rho(h)_{\nu_j}$ denotes the projection from $M$ to the $\nu_j$-eigenspace of $h$ in $M$. Throught this paper we denote by $G$ to the Hurwitz group of genus 17.

\section{The automorphism group}
In terms of generators and relations our group $G$ has the following two  inequivalent presentations (see \cite{sinkov})
\begin{align}
G&=\langle P,Q \mid P^{7}=Q^{8}=(Q P^{2})^{3}=(Q  P^{3})^{2}=Q^{4} P Q^{4} P^{4} Q^{4} P^{2}=1\rangle\label{presentationG1}\\
    &\cong\langle P,Q \mid P^{7}=Q^{8}=(Q P^{2})^{3}=(Q  P^{3})^{2}=Q^{4} P Q^{4} P^{2} Q^{4} P^{4}=1\rangle.
\end{align}
Each presentation  induces a pair of generators U,T given by
\begin{align}
 U=P^{2}Q,\\
T=P^{3}Q,
\end{align}
such that $U^{2}=T^{3}=(T^{2}U)^{7}=1$ from which one deduce that  $G$ is a Hurwitz group. 
Moreover, using an isomorphism between the two presentations one has two generating triplets 
$(U,T,T^{2}U)$, $(U',T',T'^{2}U')$ such that $T^{2}U$ and $T'^{2}U'$ belong to different  conjugacy classes. The outer automorphism group of $G$ has order 2 and
the two conjugacy classes of elements of order 7 are fixed under the action of the outer automorphisms of $G$.
This gives rise to the two non isomorphic genus 17 curves $X_{1}$ and $X_{2}$ with automorphism group $G$ (see for instance Theorem 2.17 in \cite{volklein}).\\

In this paper we have chosen the presentation (\ref{presentationG1}) as our definition of $G$.
Using Gap one can check that the following gives us a representation of $G$ as a permutation group
 \begin{align}
  P&= ( 1,13, 2,11, 4, 5, 8)( 3,10, 6,14, 7, 9,12),\\
Q&= ( 1, 7, 3, 4)( 2,11,13, 9, 6,14,10, 5),
 \end{align}
that we use to represent some elements of $G$.
Since Gap some times produces random character tables for a given group, we fix here the character table of $G$ (see Table \ref{ctblG})that we used in our calculations. 
Representatives of the conjugacy classes are shown in Table \ref{tbrepcclasses} together with the size of the normalizer group in $G$ of the cyclic group generated by a representative.
\begin{table}[!h]
\caption{Character Table of G\label{ctblG}}
{\scriptsize
\begin{displaymath}
\begin{array}{l}
\begin{array}{lrrrrrrrrrrr}
\hline\vspace{-.2cm}\\
&1A&2A&4A&2B&4B&3A&6A&8A&8B&7A&7B\\ 
\vspace{-.3cm}\\
\hline\vspace{-.2cm}\\
\chi_{1} &  1& 1& 1& 1& 1& 1& 1& 1& 1& 1& 1 \\ 
\chi_{2} &   3& 3& -1& -1& -1& 0& 0& 1& 1& \beta_{3} & \beta_{1} \\ 
\chi_{3} &    3& 3& -1& -1& -1& 0& 0& 1& 1& \beta_{1}& \beta_{3} \\ 
\chi_{4} &    6& 6& 2& 2& 2& 0& 0& 0& 0& -1& -1 \\ 
\chi_{5} &    7& 7& -1& -1& -1& 1& 1& -1& -1& 0& 0 \\ 
\chi_{6} &    7& -1& 3& -1& -1& 1& -1& 1& -1& 0& 0 \\ 
\chi_{7} &    7& -1& -1& -1& 3& 1& -1& -1& 1& 0& 0 \\ 
\chi_{8} &    8& 8& 0& 0& 0& -1& -1& 0& 0& 1& 1 \\ 
\chi_{9} &    14& -2& 2& -2& 2& -1& 1& 0& 0& 0& 0 \\ 
\chi_{10} &    21& -3& 1& 1& -3& 0& 0& -1& 1& 0& 0 \\ 
\chi_{11} &   21& -3& -3& 1& 1& 0& 0& 1& -1& 0& 0 \\
\vspace{-.3cm}\\
\hline\vspace{-.2cm}\\
\end{array}\\
\begin{array}{l}
\mbox{\scriptsize Here $\beta_{j}= \zeta^{j}+\zeta^{2j}+\zeta^{4j},~~\zeta= e^{2\pi i/7}.$}
\end{array}\\
\end{array}\\
\end{displaymath} 
}
\end{table}
\begin{table}[!h]
\caption{Representatives of the conjugacy classes of G\label{tbrepcclasses}}
{\scriptsize
\begin{displaymath}
\begin{array}{l}
\begin{array}{lcr}
\hline\vspace{-.2cm}\\
\mbox{Conjugacy class} &\mbox{ Representative} & N_{G}(<h>) \\
\vspace{-.3cm}\\
\hline\vspace{-.2cm}\\
 ~~~~1A& \mbox{Identity} & 1344  ~~~ \\
  ~~~~2A& ( 4, 7)( 5, 9)( 8,12)(10,13)& 192  ~~~ \\
  ~~~~4A& ( 4, 5, 7, 9)( 8,10,12,13)& 64   ~~~\\
  ~~~~2B& ( 2, 6)( 4, 9)( 5, 7)( 8,13)(10,12)(11,14)& 16   ~~~\\
  ~~~~4B& ( 1, 3)( 4, 5, 7, 9)( 8,13,12,10)(11,14)& 64  ~~~ \\
  ~~~~3A& ( 2, 4, 5)( 6, 7, 9)( 8,10,11)(12,13,14)& 12  ~~~ \\
  ~~~~6A& ( 1, 3)( 2, 4, 5)( 6, 7, 9)( 8,13,11,12,10,14)& 12  ~~~ \\
  ~~~~8A& ( 2, 4,14,12, 6, 7,11, 8)( 5,13, 9,10)& 32  ~~~ \\
  ~~~~8B& ( 1, 3)( 2, 4,11,12, 6, 7,14, 8)( 5,10)( 9,13)& 32   ~~~\\
  ~~~~7A& ( 1, 2, 4, 8,13,11, 5)( 3, 6, 7,12,10,14, 9)& 21   ~~~\\
  ~~~~7B& ( 1, 8, 5, 4,11, 2,13)( 3,12, 9, 7,14, 6,10)& 21  ~~~ \\
\vspace{-.3cm}\\
\hline\vspace{-.2cm}\\
\end{array}\\
\begin{array}{l}
\end{array}\\
\end{array}\\
\end{displaymath}
} 
\end{table}

\begin{table}[!h]
\caption{The characteristic polynomials of $h_{\mid_{\textsf{W}_{i}}}$\label{charpols}}
{\scriptsize
\begin{displaymath}
\begin{array}{l}
\begin{array}{cc}
\hline\vspace{-.2cm}\\
\mbox{} ~~~~~~~~~~~~~~~~~~~~~~~& h\in \mbox{7B} \\
\vspace{-.3cm}\\
\hline\vspace{-.2cm}\\
\textsf{W}_{2}&(\lambda-\zeta ) (\lambda-\zeta^2) (\lambda-\zeta^4) \\
\textsf{W}_{3}&(\lambda-\zeta^3) (\lambda-\zeta^5) (\lambda-\zeta^6) \\
\textsf{W}_{4}&(\lambda^7-1)/(\lambda-1)\\
\textsf{W}_{5},\textsf{W}_{6},\textsf{W}_{7}&(\lambda^7-1)\\
\textsf{W}_{8}&(\lambda^7-1)(\lambda-1)\\
\textsf{W}_{9}&(\lambda^7-1)^2\\
\textsf{W}_{10},\textsf{W}_{11}&(\lambda^7-1)^3\\
\vspace{-.3cm}\\
\hline\vspace{-.2cm}\\
\end{array}\\
\begin{array}{l}
\mbox{Here $\zeta=e^{2i\pi /7}$.}
\end{array}\\
\end{array}\\
\end{displaymath} 
}
\end{table}

\newpage
\section{The action on the canonical line bundle}\label{sctnx1}
One can characterize the Hurwitz curves $X_{1}$ and $X_{2}$ in terms of the action of an element of order $7$ on their canonical line bundles $K_{X_{j}}$, j=1,2. Let $d_{X_{j}}h$ be as defined in \S \ref{sectionnotation}. 
\begin{Lemma}
 Let $h\in 7B$. For  $j=1,2$ we can assume that $d_{X_{j}}h = (\zeta^{3^{j-1}}, \zeta^{3^{j-1}\cdot 2},\zeta^{3^{j-1}\cdot 4})$, where $\zeta = e^{2i\pi/7}$.
\end{Lemma}
The proof is similar to that of Lemma 2.2 in \cite{hcg14}.

So, from now on we denote by $X_{1}$ to the curve such that the action of a representative ${h_{7B}}$ of the conjugacy 
class 7B  gives the vector $d_{X_{1}}{h_{7B}}=(\zeta^1,\zeta^2,\zeta^4)$ and by $X_{2}$ to the curve 
such that  $d_{X_{2}}{h_{7B}}=(\zeta^3,\zeta^5,\zeta^6)$. Now one can use \cite[Theorem 1]{macbeath2} to complete the Table \ref{fxptscclases2dh},   and then  compute the multiplicities
$a_{i}(d)$ of the representations $\textsf{W}_{i}$, $i=1,\dots,11$, in $H^{0}(X_{j},K_{X_{j}}^{d})$ (see Tables \ref{tabla:desh0KnX1} and \ref{tabla:desh0KnX2}) by using the Atiyah--Bott fixed-point theorem and the character table of $G$.
\begin{table}[!h]
\caption{The fixed points and $d_{X_{j}}(h)$\label{fxptscclases2dh}}
{\scriptsize
\begin{displaymath}
\begin{array}{l}
\begin{array}{lcc}
\hline\vspace{-.2cm}\\
\mbox{ Conjugacy class of $h$} &\mbox{Number of fixed points in $X_{j}$} & d_{X_{j}}(h) \\
\vspace{-.3cm}\\
\hline\vspace{-.2cm}\\
  ~~~~2A& 0& --\\
  ~~~~4A& 0& --\\
  ~~~~2B& 8 & (-1,-1,-1,-1,-1,-1,-1,-1)\\
  ~~~~4B& 0&-- \\
  ~~~~3A& 4& (\omega,\omega,\omega^{2},\omega^{2}) \\
  ~~~~6A&0 &-- \\
  ~~~~8A& 0& --\\
  ~~~~8B& 0&  --\\
  ~~~~7A& 3& (\zeta^{3^{j}},\zeta^{3^{j}\cdot 2},\zeta^{3^{j}\cdot 4})\\
  ~~~~7B& 3 &(\zeta^{3^{j-1}},\zeta^{3^{j-1}\cdot 2},\zeta^{3^{j-1}\cdot 4})\\
\vspace{-.3cm}\\
\hline\vspace{-.2cm}\\
\end{array}\\
\begin{array}{l}
\end{array}\\
\mbox{Here $\omega=e^{2i\pi/3}$.}
\end{array}\\
\end{displaymath}
} 
\end{table}

\begin{table}[!h]
\caption{Decomposition of $H^{0}(X_{1},K_{X_{1}}^{d})$}\label{tabla:desh0KnX1}
{\scriptsize
\begin{displaymath}
\begin{array}{rrrrrrrrrrrrrrrrrrrrr}  
\hline\vspace{-.2cm}\\
d&1&2&3&4&5&6&7&8&9&10&11&12&13&14&15&16&17&18&19&20\\
\hline\vspace{-.2cm}\\
a_{1}(d)&   0& 0& 0& 0& 0& 1& 0& 0& 0& 0& 0& 1& 0& 1& 0& 0& 0& 1& 0& 1 \\ 
a_{2}(d)&    1& 0& 0& 0& 0& 0& 1& 1& 1& 0& 1& 0& 1& 1& 2& 1& 1& 1& 1& 1 \\ 
a_{3}(d)&    0& 0& 1& 0& 1& 0& 0& 0& 1& 1& 1& 1& 1& 0& 1& 1& 2& 1& 2& 1 \\ 
a_{4}(d)&    0& 1& 0& 1& 0& 1& 0& 2& 1& 2& 1& 2& 1& 2& 2& 3& 2& 3& 2& 3 \\ 
a_{5}(d)&    0& 0& 1& 0& 1& 1& 1& 1& 2& 1& 2& 2& 2& 2& 3& 2& 3& 3& 3& 3 \\ 
a_{6}(d)&    0& 0& 1& 0& 1& 1& 1& 1& 2& 1& 2& 2& 2& 2& 3& 2& 3& 3& 3& 3 \\ 
a_{7}(d)&    0& 0& 1& 0& 1& 1& 1& 1& 2& 1& 2& 2& 2& 2& 3& 2& 3& 3& 3& 3 \\ 
a_{8}(d)&    0& 0& 0& 1& 1& 1& 2& 1& 1& 2& 2& 2& 3& 3& 2& 3& 3& 3& 4& 4 \\ 
a_{9}(d)&   1& 0& 1& 1& 2& 1& 3& 2& 3& 3& 4& 3& 5& 4& 5& 5& 6& 5& 7& 6 \\ 
a_{10}(d)&   0& 1& 1& 2& 2& 3& 3& 4& 4& 5& 5& 6& 6& 7& 7& 8& 8& 9& 9& 10 \\ 
a_{11}(d)&    0& 1& 1& 2& 2& 3& 3& 4& 4& 5& 5& 6& 6& 7& 7& 8& 8& 9& 9& 10 \\
\end{array}
\end{displaymath}
}
\end{table}

\begin{table}[!h]
\caption{Decomposition of $H^{0}(X_{2},K_{X_{2}}^{d})$}\label{tabla:desh0KnX2}
{\scriptsize
\begin{displaymath}
\begin{array}{rrrrrrrrrrrrrrrrrrrrr}  
\hline\vspace{-.2cm}\\
d&1&2&3&4&5&6&7&8&9&10&11&12&13&14&15&16&17&18&19&20\\
\hline\vspace{-.2cm}\\
 a_{1}(d)&0& 0& 0& 0& 0& 1& 0& 0& 0& 0& 0& 1& 0& 1& 0& 0& 0& 1& 0& 1 \\ 
 a_{2}(d)&   0& 0& 1& 0& 1& 0& 0& 0& 1& 1& 1& 1& 1& 0& 1& 1& 2& 1& 2& 1 \\ 
  a_{3}(d)&  1& 0& 0& 0& 0& 0& 1& 1& 1& 0& 1& 0& 1& 1& 2& 1& 1& 1& 1& 1 \\ 
  a_{4}(d)&  0& 1& 0& 1& 0& 1& 0& 2& 1& 2& 1& 2& 1& 2& 2& 3& 2& 3& 2& 3 \\ 
  a_{5}(d)&  0& 0& 1& 0& 1& 1& 1& 1& 2& 1& 2& 2& 2& 2& 3& 2& 3& 3& 3& 3 \\ 
  a_{6}(d)&  0& 0& 1& 0& 1& 1& 1& 1& 2& 1& 2& 2& 2& 2& 3& 2& 3& 3& 3& 3 \\ 
  a_{7}(d)&  0& 0& 1& 0& 1& 1& 1& 1& 2& 1& 2& 2& 2& 2& 3& 2& 3& 3& 3& 3 \\ 
  a_{8}(d)&  0& 0& 0& 1& 1& 1& 2& 1& 1& 2& 2& 2& 3& 3& 2& 3& 3& 3& 4& 4 \\ 
 a_{9}(d)&   1& 0& 1& 1& 2& 1& 3& 2& 3& 3& 4& 3& 5& 4& 5& 5& 6& 5& 7& 6 \\ 
 a_{10}(d)&   0& 1& 1& 2& 2& 3& 3& 4& 4& 5& 5& 6& 6& 7& 7& 8& 8& 9& 9& 10 \\ 
  a_{11}(d)&  0& 1& 1& 2& 2& 3& 3& 4& 4& 5& 5& 6& 6& 7& 7& 8& 8& 9& 9& 10 \\
\end{array}
\end{displaymath}
}
\end{table}

\begin{table}[!h]
\caption{Decomposition of $S^d(\textsf{W}_{9}+\textsf{W}_{2})$}\label{tabla:desSnW9W2}
{\scriptsize
\begin{displaymath}
\begin{array}{crrrrrrrrrrr}
\hline\vspace{-.2cm}\\
d&0&1&2&3&4&5&6&7&8&9&10\\
\hline\vspace{-.2cm}\\
a_{1}(d)& 1& 0& 1& 0& 9& 5& 85& 137& 637& 1385& 4210 \\
a_{2}(d)& 0& 1& 0& 4& 6& 57& 145& 594& 1565& 4716& 11616 \\
a_{3}(d)& 0& 0& 0& 5& 6& 57& 145& 594& 1562& 4716& 11619 \\
a_{4}(d)& 0& 0& 3& 1& 34& 71& 387& 1010& 3473& 8843& 24274 \\
a_{5}(d)& 0& 0& 0& 8& 20& 118& 374& 1324& 3761& 10814& 27442 \\
a_{6}(d)& 0& 0& 0& 8& 19& 117& 367& 1317& 3743& 10796& 27391 \\
a_{7}(d)& 0& 0& 0& 8& 19& 117& 367& 1317& 3743& 10796& 27391 \\
a_{8}(d)& 0& 0& 1& 5& 32& 123& 447& 1470& 4398& 12171& 31680 \\
a_{9}(d)& 0& 1& 0& 13& 41& 234& 725& 2646& 7486& 21564& 54810 \\
a_{10}(d)& 0& 0& 3& 13& 80& 306& 1187& 3783& 11569& 31765& 83234 \\
a_{11}(d)& 0& 0& 3& 13& 80& 306& 1187& 3783& 11569& 31765& 83234 \\
\end{array}
\end{displaymath}
}
\end{table}

\begin{table}[!h]
\caption{Decomposition of $S^d(\textsf{W}_{9}+\textsf{W}_{3})$}\label{tabla:desSnW9W3}
{\scriptsize
\begin{displaymath}
\begin{array}{crrrrrrrrrrr}
\hline\vspace{-.2cm}\\
d&0&1&2&3&4&5&6&7&8&9&10\\
\hline\vspace{-.2cm}\\
a_{1}(d)&  1& 0& 1& 0& 9& 5& 85& 137& 637& 1385& 4210 \\
a_{2}(d)&  0& 0& 0& 5& 6& 57& 145& 594& 1562& 4716& 11619 \\
a_{3}(d)&  0& 1& 0& 4& 6& 57& 145& 594& 1565& 4716& 11616 \\
a_{4}(d)&  0& 0& 3& 1& 34& 71& 387& 1010& 3473& 8843& 24274 \\
a_{5}(d)&  0& 0& 0& 8& 20& 118& 374& 1324& 3761& 10814& 27442 \\
a_{6}(d)&  0& 0& 0& 8& 19& 117& 367& 1317& 3743& 10796& 27391 \\
a_{7}(d)&  0& 0& 0& 8& 19& 117& 367& 1317& 3743& 10796& 27391 \\
a_{8}(d)&  0& 0& 1& 5& 32& 123& 447& 1470& 4398& 12171& 31680 \\
a_{9}(d)&  0& 1& 0& 13& 41& 234& 725& 2646& 7486& 21564& 54810 \\
a_{10}(d)&  0& 0& 3& 13& 80& 306& 1187& 3783& 11569& 31765& 83234 \\
a_{11}(d)&  0& 0& 3& 13& 80& 306& 1187& 3783& 11569& 31765& 83234 \\
\end{array}
\end{displaymath}
}
\end{table}

\begin{table}[!h]
\caption{Decomposition of $S^d\textsf{W}_{9}$\label{decsdW9}}
{\scriptsize
\begin{displaymath}
\begin{array}{crrrrrrrrrrr}
\hline\vspace{-.2cm}\\
d&0&1&2&3&4&5&6&7&8&9&10\\
\hline\vspace{-.2cm}\\
1& 1& 0& 1& 0& 4& 2& 34& 40& 183& 332& 922 \\
2& 0& 0& 0& 2& 3& 22& 53& 184& 433& 1140& 2501 \\
3& 0& 0& 0& 2& 3& 22& 53& 184& 433& 1140& 2501 \\
4& 0& 0& 2& 0& 18& 28& 146& 308& 980& 2120& 5284 \\
5& 0& 0& 0& 4& 10& 46& 140& 408& 1047& 2612& 5922 \\
6& 0& 0& 0& 5& 9& 50& 133& 417& 1030& 2643& 5878 \\
7& 0& 0& 0& 5& 9& 50& 133& 417& 1030& 2643& 5878 \\
8& 0& 0& 1& 2& 17& 48& 165& 454& 1230& 2928& 6863 \\
9& 0& 1& 0& 7& 20& 100& 260& 843& 2060& 5267& 11772 \\
10& 0& 0& 2& 8& 39& 130& 431& 1200& 3200& 7740& 17928 \\
11& 0& 0& 2& 8& 39& 130& 431& 1200& 3200& 7740& 17928 \\
\end{array}
\end{displaymath}
}
\end{table}

\begin{table}[!h]
\caption{Decomposition of $S^d\textsf{W}_{2}$\label{decsdW2}}
{\scriptsize
\begin{displaymath}
\begin{array}{crrrrrrrrrrr}
\hline\vspace{-.2cm}\\
d&0&1&2&3&4&5&6&7&8&9&10\\
\hline\vspace{-.2cm}\\
a_{1}(d)& 1& 0& 0& 0& 1& 0& 1& 0& 1& 0& 1 \\
a_{2}(d)& 0& 1& 0& 0& 0& 1& 0& 1& 1& 2& 0 \\
a_{3}(d)& 0& 0& 0& 1& 0& 1& 0& 1& 0& 2& 1 \\
a_{4}(d)& 0& 0& 1& 0& 1& 0& 2& 0& 3& 1& 4 \\
a_{5}(d)& 0& 0& 0& 1& 0& 1& 1& 2& 1& 3& 2 \\
a_{6}(d)& 0& 0& 0& 0& 0& 0& 0& 0& 0& 0& 0 \\
a_{7}(d)& 0& 0& 0& 0& 0& 0& 0& 0& 0& 0& 0 \\
a_{8}(d)& 0& 0& 0& 0& 1& 1& 1& 2& 2& 2& 3 \\
a_{9}(d)& 0& 0& 0& 0& 0& 0& 0& 0& 0& 0& 0 \\
a_{10}(d)& 0& 0& 0& 0& 0& 0& 0& 0& 0& 0& 0 \\
a_{11}(d)& 0& 0& 0& 0& 0& 0& 0& 0& 0& 0& 0 \\
\end{array}
\end{displaymath}
}
\end{table}

\begin{table}[!h]
\caption{Decomposition of $S^d\textsf{W}_{3}$\label{decsdW3}}
{\scriptsize
\begin{displaymath}
\begin{array}{crrrrrrrrrrr}
\hline\vspace{-.2cm}\\
d&0&1&2&3&4&5&6&7&8&9&10\\
\hline\vspace{-.2cm}\\
a_{1}(d)& 1& 0& 0& 0& 1& 0& 1& 0& 1& 0& 1 \\
a_{2}(d)& 0& 0& 0& 1& 0& 1& 0& 1& 0& 2& 1 \\
a_{3}(d)& 0& 1& 0& 0& 0& 1& 0& 1& 1& 2& 0 \\
a_{4}(d)& 0& 0& 1& 0& 1& 0& 2& 0& 3& 1& 4 \\
a_{5}(d)& 0& 0& 0& 1& 0& 1& 1& 2& 1& 3& 2 \\
a_{6}(d)& 0& 0& 0& 0& 0& 0& 0& 0& 0& 0& 0 \\
a_{7}(d)& 0& 0& 0& 0& 0& 0& 0& 0& 0& 0& 0 \\
a_{8}(d)& 0& 0& 0& 0& 1& 1& 1& 2& 2& 2& 3 \\
a_{9}(d)& 0& 0& 0& 0& 0& 0& 0& 0& 0& 0& 0 \\
a_{10}(d)& 0& 0& 0& 0& 0& 0& 0& 0& 0& 0& 0 \\
a_{11}(d)& 0& 0& 0& 0& 0& 0& 0& 0& 0& 0& 0 \\
\end{array}
\end{displaymath}
}
\end{table}

\newpage
  \section{The quadratic equations and the matrix representation of G}\label{sectionquadrics}
From Tables \ref{tabla:desh0KnX1} and \ref{tabla:desh0KnX2} we have
\begin{align}
 H^0(X_{i},K_{X_{i}})\:&=\:\:\textsf{W}_{9}\oplus \textsf{W}_{i+1}, \:\:\mbox{ i=1,2, }\\
&\mbox{ and }\nonumber\\
H^0(X_{i},K_{X_{i}}^{\otimes2})\:&=\:\:\:\:\:\:\textsf{W}_{4}\oplus \textsf{W}_{10}\oplus \textsf{W}_{11}.
\end{align}
So we assume that $X_{i}\subset\mathbb{P}^{16}= \mathbb{P}(\textsf{V}_{i}^{*})$,
 where $\textsf{V}_{i}$ has the structure of the representation
$\textsf{W}_{9}\oplus \textsf{W}_{i+1}$. From Tables  \ref{tabla:desSnW9W2} and \ref{tabla:desSnW9W3} we have
\begin{align}
S^2(\textsf{V}_{i} )\:&=\:\textsf{W}_{1}\oplus \textsf{W}_{4}^{3}\oplus \textsf{W}_{8}\oplus \textsf{W}_{10}^{3}\oplus \textsf{W}_{11}^3,\mbox{ i=1,2}.
\end{align}
It follows that
\begin{equation}
  H^0(\mathbb{P}^{16},\mathcal{I}_{X_{i}}(2))=\textsf{W}_{1}\oplus \textsf{W}_{4}^2\oplus \textsf{W}_{8}\oplus \textsf{W}_{10}^2\oplus \textsf{W}_{11}^2,\mbox{ i=1,2}.
\end{equation}
Where $\mathcal{I}_{X_{i}}$ denotes the ideal sheaf of $X_{i}$ in $\mathbb{P}^{16}$.
From now on we will use the notation
\[ I_{X_{i}}(2)=H^0(\mathbb{P}^{16},\mathcal{I}_{X_{i}}(2)).\]
We have 
\begin{align}
\pi_{{}_{\textsf{W}_{1}}}(I_{X_{i}}(2))= \pi_{{}_{\textsf{W}_{1}}}(S^{2}\textsf{V}_{i}),\\
\mbox{ and ~~~~~~~~~~~~~~~}\nonumber\\
\pi_{{}_{\textsf{W}_{8}}}(I_{X_{i}}(2))= \pi_{{}_{\textsf{W}_{8}}}(S^{2}\textsf{V}_{i}).
 \end{align}
The remaining spaces $\pi_{{}_{\textsf{W}_{4}}}(I_{X_{i}}(2))$,
 $\pi_{{}_{\textsf{W}_{10}}}(I_{X_{i}}(2))$ and $\pi_{{}_{\textsf{W}_{11}}}(I_{X_{i}}(2))$
will be found in Sections \ref{secpartW4},\ref{secpartW10} and \ref{secpartW11} below.\\
For calculations we fix an action of $G$ on $\mathbb{P}^{16}$. 
 In $\textsf{V}_{i}$ (respectively $\textsf{V}_{i}^{*}$ ) we consider a basis 
$\{y_{1},\dots,y_{17}\}$ (respectively the dual basis $\{e_{1},\dots,e_{17}\}$). 
Sometimes we   use the following subspaces of $V_{i}$
\begin{align}
\textsf{V}_{{}_{i,9}}\:\: &=\:\:<y_{1},\dots,y_{14}>,\\
\textsf{V}_{{}_{i,i+1}}& =\:<y_{15},y_{16},y_{17}>.
\end{align}

With respect to these bases, the action of the generators $P$ and $Q$ on 
 $\textsf{V}_{i}$ (respectively $\textsf{V}_{i}^{*}$) is  given by the matrices
$\mathrm{diag}(P_{{}_{\textsf{W}_{9}}}, P_{{}_{\textsf{W}_{i+1}}})$ and 
$\mathrm{diag}(Q_{{}_{\textsf{W}_{9}}}, Q_{{}_{\textsf{W}_{i+1}}})$ 
(respectively the inverse  transpose of these matrices), 
where the  matrices $P_{{}_{\textsf{W}_{j}}}$, $Q_{{}_{\textsf{W}_{j}}}$, $j=9,2$ are 
defined in formulae (\ref{mpx2}), (\ref{mqx2}), (\ref{mpx9}) and(\ref{mqx9}) (recall that $\omega=e^{2i\pi/3}$ and $\zeta=e^{2i\pi/7}$). 
The matrices $P_{{}_{\textsf{W}_{3}}}$, $Q_{{}_{\textsf{W}_{3}}}$ 
could be taken as the respective inverse transpose of $P_{{}_{\textsf{W}_{2}}}$, $Q_{{}_{\textsf{W}_{2}}}$ 
however we will focus only on the curve $X_{1}$.\\
\begin{Remark}\label{rmk2}
 \textnormal{ One can identify $G$ with a subgroup of $Aut(G)$. The action of $G$ extends to an action of $Aut(G)$ on $\textsf{V}_{1}$. So the exterior automorphisms move $X_{1}\subset \mathbb{P}(\textsf{V}_{1}^{*})$ giving rise to another copy  $X_{1}'\subset \mathbb{P}(\textsf{V}_{1}^{*})$ of $X_{1}$. For instance, we can assume $Aut(G)=< P,Q,E>$, where $E$ is such that the action on $\textsf{V}_{1}$ is given by 
the matrix $\mathrm{diag}(E_{\textsf{W}_{9}},E_{\textsf{W}_{2}})$, where $E_{\textsf{W}_{9}}$ and $E_{\textsf{W}_{2}}$ are defined in formulae
(\ref{sw9}) and (\ref{sw2}) respectively.
}
\end{Remark}

\begin{Remark}\label{rmk1}
 \textnormal{
Let $Y$ denote the Klein quartic. We have that $Y$ is isomorphic to the quotient $X_{i}/N$ for $i=1,2$, where $N$ is the normal subgroup of $G$ generated by the involutions in the class $2A$. Let $f_{i}:X_{i}\mapsto Y$ be the quotient maps. Then we have that  the action of $G$ induced on $Y$ from $X_{1}$ or $X_{2}$ is such that ${d_{Y}}{h_{7B}}=(\zeta^1,\zeta^2,\zeta^4)$ or ${d_{Y}}{h_{7B}}=(\zeta^3,\zeta^5,\zeta^6)$ respectively. The action induced  by $f_{1}$ or $f_{2}$ gives $H^{0}(Y,K_{Y})$ the structure of  the representation $\textsf{W}_{2}$ or $\textsf{W}_{3}$ respectively.
Since $f_{i}^{*}K_{Y}\cong K_{X_{i}}$,
 we see that  $f_{i}$ can be recovered by composing the canonical embedding
   $X_{i}\subset\mathbb{P}(\textsf{V}_{i}^{*})$
    with  the projection 
    $\pi_{i}:\mathbb{P}(\textsf{V}_{i}^{*})\dashrightarrow \mathbb{P}( \textsf{V}_{{}_{i,i+1}}^{*})$.
}
\end{Remark}

{\scriptsize
\begin{align}
P_{{}_{\textsf{W}_{2}}}&\,=\:\qquad\qquad\qquad\left(\begin{array}{ccc}
  \zeta^3& 0& 0 \\  0& \zeta^5& 0 \\  0& 0& \zeta^6 \\
\end{array}\right)\label{mpx2}\\
\nonumber\\
Q_{{}_{\textsf{W}_{2}}}&\,=\:\frac{1}{\sqrt{-7}}\left(\begin{array}{ccc}
  \zeta^4-\zeta^6& -1+\zeta^3& \zeta-\zeta^2 \\ 
   -\zeta+\zeta^4& \zeta^2-\zeta^3& -1+\zeta^5\\ 
   -1+\zeta^6& \zeta^2-\zeta^4& \zeta-\zeta^5 \\
\end{array}\right)\label{mqx2}
\end{align}
}

{\scriptsize
\begin{equation}\label{mpx9}
P_{{}_{\textsf{W}_{9}}}=\left(\begin{array}{cccccccccccccc}
  0& 0& 0& 0& 0& 0& 0& 0& 0& 0& 0& 0& 0& \omega^2 \\  0& 0& 0& 0& 0& 0& 0& 0& 0& 0& 0& 0& \omega& 0 \\  0& 0& 0& 0& 0& 0& 1& 0& 0& 0& 0& 0& 0& 0 \\  0& 0& 0& 0& 0& 0& 0& 1& 0& 0& 0& 0& 0& 0 \\  0& 0& 0& 0& 0& 0& 0& 0& -\omega^2& 0& 0& 0& 0& 0 \\  0& 0& 0& 0& 0& 0& 0& 0& 0& -\omega& 0& 0& 0& 0 \\  0& 0& 0& 0& -\omega^2& 0& 0& 0& 0& 0& 0& 0& 0& 0 \\  0& 0& 0& 0& 0& -\omega& 0& 0& 0& 0& 0& 0& 0& 0 \\  0& -\omega^2& 0& 0& 0& 0& 0& 0& 0& 0& 0& 0& 0& 0 \\  -\omega& 0& 0& 0& 0& 0& 0& 0& 0& 0& 0& 0& 0& 0 \\  0& 0& \omega& 0& 0& 0& 0& 0& 0& 0& 0& 0& 0& 0 \\  0& 0& 0& \omega^2& 0& 0& 0& 0& 0& 0& 0& 0& 0& 0 \\  0& 0& 0& 0& 0& 0& 0& 0& 0& 0& -\omega& 0& 0& 0 \\  0& 0& 0& 0& 0& 0& 0& 0& 0& 0& 0& -\omega^2& 0& 0 
\end{array}\right)
\end{equation}
}
{\scriptsize
\begin{equation}\label{mqx9}
Q_{{}_{\textsf{W}_{9}}}=\left(\begin{array}{cccccccccccccc}
 0& 0& 0& 0& 0& 0& 0& 1& 0& 0& 0& 0& 0& 0 \\  0& 0& 0& 0& 0& 0& 1& 0& 0& 0& 0& 0& 0& 0 \\  0& 0& 0& 0& 0& 0& 0& 0& 1& 0& 0& 0& 0& 0 \\  0& 0& 0& 0& 0& 0& 0& 0& 0& 1& 0& 0& 0& 0 \\  0& 0& 0& 0& 0& \omega^2& 0& 0& 0& 0& 0& 0& 0& 0 \\  0& 0& 0& 0& \omega& 0& 0& 0& 0& 0& 0& 0& 0& 0 \\  0& 0& 0& 0& 0& 0& 0& 0& 0& 0& 0& 0& 0& -\omega \\  0& 0& 0& 0& 0& 0& 0& 0& 0& 0& 0& 0& -\omega^2& 0 \\  0& 0& 0& -\omega^2& 0& 0& 0& 0& 0& 0& 0& 0& 0& 0 \\  0& 0& -\omega& 0& 0& 0& 0& 0& 0& 0& 0& 0& 0& 0 \\  0& -\omega^2& 0& 0& 0& 0& 0& 0& 0& 0& 0& 0& 0& 0 \\  -\omega& 0& 0& 0& 0& 0& 0& 0& 0& 0& 0& 0& 0& 0 \\  0& 0& 0& 0& 0& 0& 0& 0& 0& 0& 0& -1& 0& 0 \\  0& 0& 0& 0& 0& 0& 0& 0& 0& 0& -1& 0& 0& 0 
\end{array}\right)
\end{equation}
}

{\scriptsize
\begin{align}\label{sw9}
E_{{}_{\textsf{W}_{9}}}=\left(\begin{array}{cccccccccccccc}
  0& 0& 0& 0& 0& 0& 1& 0& 0& 0& 0& 0& 0& 0 \\ 
   0& 0& 0& 0& 0& 0& 0& -1& 0& 0& 0& 0& 0& 0 \\ 
   0& 0& 0& \omega^2& 0& 0& 0& 0& 0& 0& 0& 0& 0& 0 \\ 
   0& 0& -\omega& 0& 0& 0& 0& 0& 0& 0& 0& 0& 0& 0 \\ 
   0& 0& 0& 0& 0& -\omega^2& 0& 0& 0& 0& 0& 0& 0& 0 \\ 
   0& 0& 0& 0& \omega& 0& 0& 0& 0& 0& 0& 0& 0& 0 \\ 
   -1& 0& 0& 0& 0& 0& 0& 0& 0& 0& 0& 0& 0& 0 \\ 
   0& 1& 0& 0& 0& 0& 0& 0& 0& 0& 0& 0& 0& 0 \\ 
   0& 0& 0& 0& 0& 0& 0& 0& -1& 0& 0& 0& 0& 0 \\ 
   0& 0& 0& 0& 0& 0& 0& 0& 0& 1& 0& 0& 0& 0 \\ 
   0& 0& 0& 0& 0& 0& 0& 0& 0& 0& 0& 0& -\omega& 0 \\ 
   0& 0& 0& 0& 0& 0& 0& 0& 0& 0& 0& 0& 0& \omega^2 \\ 
   0& 0& 0& 0& 0& 0& 0& 0& 0& 0& -\omega^2& 0& 0& 0 \\ 
   0& 0& 0& 0& 0& 0& 0& 0& 0& 0& 0& \omega& 0& 0 \\
\end{array}\right)\\
\end{align}
}

{\scriptsize
\begin{align}\label{sw2}
E_{{}_{\textsf{W}_{2}}}=\frac{1}{\sqrt{-7}}\left(\begin{array}{ccc}
  \zeta-\zeta^6& -\zeta^2+\zeta^6 &  -\zeta^2+\zeta^3  \\ 
  -\zeta+\zeta^5 & -\zeta^3+\zeta^4&  -\zeta+\zeta^3 \\ 
  -\zeta^4+\zeta^5 & -\zeta^4+\zeta^6 & \zeta^2-\zeta^5 \\
\end{array}\right)
\end{align}
}

 \section{The generators of $\pi_{{}_{\textsf{W}_{4}}}(I_{X_{1}}(2))$ }\label{secpartW4}
Let $h_{7B}=P^2\! Q^{}\! P^3\! Q =(1,8,5,4,11,2,13)(3,12,9,7,14,6,10) $ 
be the representative of the conjugacy class $7B$ of $G$. 
From Table \ref{charpols} we see that  $\rho(h_{7B})_{\zeta^2}(\textsf{W}_{4})$ has dimension $1$. 
So, since $\textsf{W}_{4}$ is an irreducible $\mathbb{C}G$-module, the subspace $\rho(h_{7B})_{\zeta^2}(\pi_{{}_{\textsf{W}_{4}}}(I_{X_{1}}(2)))$ 
generates $\pi_{{}_{\textsf{W}_{4}}}(I_{X_{1}}(2))$ as a $\mathbb{C}G$-module. 
 We shall find a point $p\in X_{1}$ such that not all  $\zeta^2$-eigenquadrics of $h_{7B}$ in 
$\pi_{{}_{\textsf{W}_{4}}}(S^2\textsf{V}_{1})$  vanish on $p$.   
Then the subspace $\rho(h_{7B})_{\zeta^2}(\pi_{{}_{\textsf{W}_{4}}}(I_{X_{1}}(2)))$ is the set of 
 all the  $\zeta^2$-eigenquadrics of $h_{7B}$ in 
$\pi_{{}_{\textsf{W}_{4}}}(S^2\textsf{V}_{1})$ that vanish on $p$.

If $p$ is a fixed point of $h_{7B}$ in $X_{1}$, then it lives in a subspace \\
$\mathbb{P}(\rho(h_{7B})_{\zeta^i}(\textsf{V}_{1}^{*}))\subset\mathbb{P}(\textsf{V}_{1}^{*})$  for some $i$. 
It can be checked that $i=6,5$ or $3$ if ${d{h_{7B}}}_{p}= \zeta^1,\zeta^2 $ or $\zeta^4$ respectively.\\
Now we will look for the fixed point $p$ of $h_{7B}$ in $X_{1}$ lying on $\mathbb{P}(\rho(h_{7B})_{\zeta^6}(\textsf{V}_{1}^{*}))$.
 Note that $\rho(h_{7B})_{\zeta^6}(\textsf{V}_{{}_{1,9}}^{*}) $ has dimension 2 (see Table \ref{charpols})
and that 
$f_{1}(p)$ (see Remark \ref{rmk1}) 
lies on the subspace  $\rho(h_{7B})_{\zeta^6}(\textsf{V}_{{}_{1,2}}^{*}) $ which is one dimensional. Therefore we can assume that 
\begin{equation}
p=\alpha\cdot v_{1}+\beta\cdot v_{2}+ v_{3},
\end{equation} where $\alpha$,$\beta\in \mathbb{C}$ and $\{v_{1},v_{2}\}$ and $\{v_{3}\}$ are bases for 
$\rho(h_{7B})_{\zeta^6}(\textsf{V}_{{}_{1,9}}^{*}) $ and
 \\
 $\rho(h_{7B})_{\zeta^6}(\textsf{V}_{{}_{1,2}}^{*}) $ respectively.

To find $\alpha$ and $\beta $ we  will evaluate on the point $p$ the $\zeta^2$-eigenquadric  
$q_{\zeta^2}$ of $h_{7B}$ in $\pi_{\textsf{W}_{8}}(I_{X_{1}}(2))$ ; namely
\begin{equation}\label{qz2x8}
 q_{\zeta^2}=\rho(h_{7B})_{\zeta^2}(\pi_{{}_{\textsf{W}_{8}}}(y_{1}^2)).
\end{equation}

We take
\begin{equation}\begin{array}{lcl}
v_{i}&=&\rho(h_{7B})_{\zeta^6}(e_{i}), ~i=1,2  \\
&&\mbox{~and}\\
v_{3}&=&\rho(h_{7B})_{\zeta^6}(e_{17}).
\end{array}
\end{equation}
 We see that 
\begin{equation} q_{\zeta^2}(v_{i})\not= 0, ~i=1,2. 
\end{equation}
We assume that $\beta=1$ since we always can move the fixed point $p$ by a $G$-equivariant change of
 coordinates $t_{\lambda}:e_{i}\rightarrow\lambda e_{i}$, $i=1,\dots,14$ and $e_{i}\rightarrow e_{i}$, $i=15,16,17$.
So we have
\begin{equation}\label{eqofxpnt}\mbox{\scriptsize
 $\alpha^{2}+1/7(-6 \xi-4 \xi^{2}+4 \xi^{4}-3 \xi^{5}-2 \xi^{8}+3 \xi^{10}+6 \xi^{11}-\xi^{13}+2 \xi^{16}+\xi^{17}-2 \xi^{19}+2 \xi^{20})=0,$
}\end{equation}

where $\xi=e^{2i\pi/21}.$\\
Notice that the corresponding points to each root of (\ref{eqofxpnt}) live in projectively equivalent curves in $\mathbb{P}(\textsf{V}_{1}^{*})$. One is obtained from the other by applying the exterior automorphism in  Remark \ref{rmk2} followed by the equivariant translation $t_{\lambda}$ with suitable
$\lambda$.

Now, a basis for the $\zeta^2$-eigenspace of ${h_{7B}}_{\mid_{\pi_{{}_{\textsf{W}_{4}}}(S^{2}\textsf{V}_{1})}}$ is given by $\{b_{1}, b_{2},b_{3}\}$, where
\begin{equation} \begin{array}{l}
                  ḅ_{1}=\rho(h_{7B})_{\zeta^2}(\pi_{{}_{\textsf{W}_{4}}}(y_{1}^2)),\\
ḅ_{2}=\rho(h_{7B})_{\zeta^2}(\pi_{{}_{\textsf{W}_{4}}}(y_{1}y_{2})),\\
ḅ_{3}=y_{17}^{2}.\\
                 \end{array}
\end{equation}
Evaluating on $p$ one proves that 
$b_{i}(p)\not=0$, $i=1,2,3.$ Thus the quadrics
\begin{equation}\label{vathetaeqs}
 \vartheta_{1}= b_{1}(p)b_{2}-b_{2}(p)b_{1} \mbox{ and } \vartheta_{2}= b_{1}(p)b_{3}-b_{3}(p)b_{1}
\end{equation}
 form a basis for the $\zeta^2$-eigenspace of ${h_{7B}}_{\mid_{\pi_{{}_{\textsf{W}_{4}}}(I_{X_{1}}(2))}}$.
\section{The generators of $\pi_{{}_{\textsf{W}_{10}}}(I_{X_{1}}(2))$}\label{secpartW10}
We shall find  linearly independent elements $\theta_{1},\theta_{2},\theta_{3}\in \pi_{{}_{\textsf{W}_{10}}}(S^2\textsf{V}_{1})$ such that
the $G$-module  generated by each one of them is isomorphic to $\textsf{W}_{10}$ and such that $\theta_{1},\theta_{2},\theta_{3}$
 are $G$-\emph{equivalent}, meaning that if $\phi_{ij}$ is an isomorphism of $G$-modules from
 the $G$-module generated by $\theta_{i}$ to that generated by $\theta_{j}$ then
 $\phi_{ij}(\thetạ_{i})=\lambda_{ij}\theta_{j}$ for some $\lambda_{ij}\in \mathbb{C}\backslash\{0\}$.
We take $\theta_{i}$ to be a $\zeta^2$-eigenquadric of ${h_{7B}}_{\mid_{\pi_{{}_{\textsf{W}_{10}}}(S^2\textsf{V}_{1})}}$. Let $p$ be the fixed point of  $h_{7B}$
 we found in Section \ref{secpartW4}. If $\theta_{i}(p)\not=0$, $i=1,2,3$ then
\begin{equation}\label{prtx10geneqns}
\begin{array}{l}
 \wp_{1}=\theta_{1}(p)\theta_{2}-\theta_{2}(p)\theta_{1},\\
 \wp_{2}=\theta_{1}(p)\theta_{3}-\theta_{3}(p)\theta_{1}
\end{array}
\end{equation}
generate $\pi_{{}_{\textsf{W}_{10}}}(I_{X_{1}}(2))$.\\
We first notice that
\begin{equation}
\begin{array}{ccl}
S^2 \textsf{V}_{1}&= &S^{2}\textsf{V}_{{}_{1,9}}\oplus S^2\textsf{V}_{{}_{1,2}}\oplus\:\textsf{V}_{{}_{1,9}}\!\!\otimes\textsf{V}_{{}_{1,2}}.\\
\end{array}
\end{equation}
 Using that (see Tables \ref{decsdW9},\ref{decsdW2})
\begin{equation}
\begin{array}{lccl}
& S^{2} \textsf{V}_{{}_{1,9}}&=&\textsf{W}_{1}\oplus\textsf{W}_{4}^2\oplus\textsf{W}_{8}\oplus\textsf{W}_{10}^2\oplus \textsf{W}_{11}^2, \\
& S^{2} \textsf{V}_{{}_{1,2}}&=&\textsf{W}_{4},\\
\end{array}
\end{equation}
one concludes that
\begin{equation}
\begin{array}{lccl}
&\textsf{V}_{{}_{1,9}}\!\!\otimes\textsf{V}_{{}_{1,2}}&=&\textsf{W}_{10}\oplus \textsf{W}_{11}.\\
\end{array}
\end{equation}
 Then we define an isomorphism of $G$-modules
\[ \varrho_{{}_{\textsf{W}_{10}}}:\pi_{{}_{\textsf{W}_{10}}}(\textsf{V}_{{}_{1,9}}\!\!\otimes\textsf{V}_{{}_{1,2}})\oplus \pi_{{}_{\textsf{W}_{10}}}(\textsf{V}_{{}_{1,9}}\!\!\otimes\textsf{V}_{{}_{1,2}})\rightarrow \pi_{{}_{\textsf{W}_{10}}}(S^{2} \textsf{V}_{{}_{1,9}}).\]
So given \[\theta_{1}\in\pi_{{}_{\textsf{W}_{10}}}(\textsf{V}_{{}_{1,9}}\!\!\otimes\textsf{V}_{{}_{1,2}})\]
 one can take
\begin{equation}
\begin{array}{ccl}
\theta_{2}&=&\varrho_{{}_{\textsf{W}_{10}}}((\theta_{1},0))\\
\theta_{3}&=&\varrho_{{}_{\textsf{W}_{10}}}((0,\theta_{1})).\\
\end{array}
\end{equation}
To construct $\varrho_{{}_{\textsf{W}_{10}}}$  we consider the bases 
$B_{1}$ of $\pi_{{}_{\textsf{W}_{10}}}(\textsf{V}_{{}_{1,9}}\!\!\otimes\textsf{V}_{{}_{1,2}})$
and $B_{2}$ of $\pi_{{}_{\textsf{W}_{10}}}(S^{2}\textsf{V}_{{}_{1,9}})$
 given in formulae (\ref{basB1}) and (\ref{basB2}).
Let $Ạ,C$ be the matrix representations of $G$ corresponding to the action on 
$\pi_{{}_{\textsf{W}_{10}}}(\textsf{V}_{{}_{1,9}}\!\!\otimes\textsf{V}_{{}_{1,2}})$
and on
$\pi_{{}_{\textsf{W}_{10}}}(S^{2}\textsf{V}_{{}_{1,9}})$ with respect to these bases. If $M_{0}$ is a $42\times42$ complex matrix 
such that $\det{S}\not=0$, where $S$ given by equation  (\ref{cnjmtrx}), then one can take  $\varrho_{{}_{\textsf{W}_{10}}}$ 
to be the linear map induced by $S$  with respect to the bases $B_{1}\dot{\cup}B_{1}$ and $B_{2}$.
\begin{equation}\label{cnjmtrx} 
S= \sum_{h\in G} C_{h}\cdot M_{0}\cdot\diag(A_{h},A_{h})^{-1}.
\end{equation}
For our calculation we use $M_{0}=I_{42}$ the identity matrix.

\begin{equation}\label{basB1}
\begin{array}{ccl}
 B_{1}&=&\{ \mbox{\small\ss{}}_{1,j},\mbox{\small\ss{}}_{2,j},\mbox{\small\ss{}}_{3,j}\mid j=1\dots 7\},\\
\end{array}\\
\end{equation}
where
\begin{equation}
\begin{array}{ccl}
\mbox{\ss{}}_{i,j}&=& \rho(h_{7B})_{\zeta^j}(Q\cdot\mbox{\ss{}}_{i,2}),\mbox{ for $j\not= 2$},\\ 

\mbox{\ss{}}_{1,2}&=&\rho(h_{7B})_{\zeta^2}(\pi_{{}_{\textsf{W}_{10}}}(y_{1}y_{15})),\\
\mbox{\ss{}}_{i+1,2}&=& \rho(h_{7B})_{\zeta^2}(Q\cdot\mbox{\ss{}}_{i,2}).\\ 

\end{array}
\end{equation}
\begin{equation}
\begin{array}{ccl}\label{basB2}
B_{2}&=&\{\kappa_{1,j},\kappa_{2,j},\kappa_{3,j},\kappa_{4,j},\kappa_{5,j},\kappa_{6,j}\mid j=1\dots 7\},
\end{array}\\
\end{equation}
where
\begin{equation}
\begin{array}{ccl}
\kappa_{i,j}&=&\rho(h_{7B})_{\zeta^j}(Q\cdot\kappa_{i,2}),\mbox{ for $j\not= 2$},\\

\kappa_{1,2}&=&\rho(h_{7B})_{\zeta^2}(\pi_{{}_{\textsf{W}_{10}}}(y_{1}y_{3})),\\
\kappa_{4,2}&=&\rho(h_{7B})_{\zeta^2}(Q\cdot\rho(h_{7B})_{\zeta^3}(Q\cdot\kappa_{1,2})),\\
\kappa_{i+1,2}&=&\rho(h_{7B})_{\zeta^2}(Q\cdot\kappa_{i,2}), \mbox{ for $i= 1,2,4,5$}.\\

\end{array}
\end{equation}

We take $\theta_{1}=\mbox{\ss{}}_{1,2}$.
Then \begin{equation}
\begin{array}{ccl}
\theta_{2}&=&\sum_{i=1}^{6}s_{i+6,4}\kappa_{i,2},\\
\theta_{3}&=&\sum_{i=1}^{6}s_{i+6,25}\kappa_{i,2}.\\
\end{array}
\end{equation}

\section{The generators of $\pi_{{}_{\textsf{W}_{11}}}(I_{X_{1}}(2))$}\label{secpartW11}
 The situation is similar to that of $\pi_{{}_{\textsf{W}_{10}}}(I_{X_{1}}(2))$. We define an isomorphism of $G$-modules

\[ \varrho_{{}_{\textsf{W}_{11}}}:\pi_{{}_{\textsf{W}_{11}}}(\textsf{V}_{{}_{1,9}}\!\!\otimes\textsf{V}_{{}_{1,2}})\oplus \pi_{{}_{\textsf{W}_{11}}}(\textsf{V}_{{}_{1,9}}\!\!\otimes\textsf{V}_{{}_{1,2}})\rightarrow \pi_{{}_{\textsf{W}_{11}}}(S^{2} \textsf{V}_{{}_{1,9}}).\]

which corresponds to a matrix
\begin{equation}\label{cnjmtrx2} 
S'= \sum_{h\in G} C_{h}'\cdot M_{0}'\cdot\diag(A_{h}',A_{h}')^{-1}.
\end{equation}
where $Ạ',C'$ are the matrix representations of $G$ corresponding to the action on 
$\pi_{{}_{\textsf{W}_{11}}}(\textsf{V}_{{}_{1,9}}\!\!\otimes\textsf{V}_{{}_{1,2}})$
and on
$\pi_{{}_{\textsf{W}_{11}}}(S^{2} \textsf{V}_{{}_{1,9}})$
with respect to fixed  bases $B_{1}'$ and $B_{2}'$ respectively( see formulae (\ref{basBp1}), (\ref{basBp2})).
We use $M_{0}'=I_{42}$ the identity matrix. 
Then $\wp_{1}'$ and $\wp_{2}'$  in formula (\ref{prtx11geneqns}) generate $\pi_{{}_{\textsf{W}_{11}}}(I_{X_{1}}(2))$.
\begin{equation}
\begin{array}{ccl}\label{basBp1}
 B_{1}'&=&\{ \mbox{\small\ss{}}_{1,j}',\mbox{\small\ss{}}_{2,j}',\mbox{\small\ss{}}_{3,j}'\mid j=1\dots 7\},\\
\end{array}\\
\end{equation}
where
\begin{equation}
\begin{array}{ccl}
\mbox{\ss{}}_{i,j}'&=& \rho(h_{7B})_{\zeta^j}(Q\cdot\mbox{\ss{}}_{i,2}'),\mbox{ for $j\not= 2$},\\ 

\mbox{\ss{}}_{1,2}'&=&\rho(h_{7B})_{\zeta^2}(\pi_{{}_{\textsf{W}_{11}}}(y_{1}y_{15})),\\
\mbox{\ss{}}_{i+1,2}'&=& \rho(h_{7B})_{\zeta^2}(Q\cdot\mbox{\ss{}}_{i,2}').\\ 

\end{array}
\end{equation}
\begin{equation}
\begin{array}{ccl}\label{basBp2}
B_{2}'&=&\{\kappa_{1,j}',\kappa_{2,j}',\kappa_{3,j}',\kappa_{4,j}',\kappa_{5,j}',\kappa_{6,j}'\mid j=1\dots 7\},
\end{array}\\
\end{equation}
where
\begin{equation}
\begin{array}{ccl}
\kappa_{i,j}'&=&\rho(h_{7B})_{\zeta^j}(Q\cdot\kappa_{i,2}'),\mbox{ for $j\not= 2$},\\
\kappa_{1,2}'&=&\rho(h_{7B})_{\zeta^2}(\pi_{{}_{\textsf{W}_{11}}}(y_{1}y_{3})),\\
\kappa_{4,2}'&=&\rho(h_{7B})_{\zeta^2}(Q\cdot\rho(h_{7B})_{\zeta^3}(Q\cdot\kappa_{1,2}')),\\
\kappa_{i+1,2}'&=&\rho(h_{7B})_{\zeta^2}(Q\cdot\kappa_{i,2}')\mbox{ for $i=1,2,4,5$}.\\
\end{array}
\end{equation}
Taking $\theta_{1}' = \mbox{\ss{}}_{1,2}'$ we take
\begin{equation}
\begin{array}{ccl}
\theta_{2}'&=&\sum_{i=1}^{6}s_{i+6,4}'\kappa_{i,2}',\\
\theta_{3}'&=&\sum_{i=1}^{6}s_{i+6,25}'\kappa_{i,2}'.\\
\end{array}
\end{equation}
One checks that $\theta_{1}'(p)\theta_{2}'(p)\theta_{3}'(p)\not=0$. So we define
\begin{equation}\label{prtx11geneqns}
\begin{array}{l}
 \wp_{1}'=\theta_{1}'(p)\theta_{2}'-\theta_{2}'(p)\theta_{1}',\\
 \wp_{2}'=\theta_{1}'(p)\theta_{3}'-\theta_{3}'(p)\theta_{1}'.
\end{array}
\end{equation}
Although one can speed up the calculation of (\ref{cnjmtrx}) and (\ref{cnjmtrx2}) by using representatives of cosets of a subgroup of G, the computation was still very slow.
\section{A basis for $I_{X_{1}}(2)$ }\label{basI2}
Now we can give formulae for the elemens of a basis for the space of quadrics in the canonical ideal of $X_{1}$.
\begin{theorem}\label{tbasI2}
The union of the sets of quadrics $D_{i}$, $i=1,4,8,10,11$ below form a basis for $H^0(\mathbb{P}^{16},\mathcal{I}_{X_{i}}(2))$.
\end{theorem}

 Let $\alpha$ be a solution of (\ref{eqofxpnt}) and $p$ be the corresponding point. Let $q_{\zeta^2}$ be as defined in (\ref{qz2x8}).
Let $\vartheta_{1}$,$\vartheta_{2}$ be as defined in (\ref{vathetaeqs}). 
Let  $\wp_{1}$,$\wp_{2}$,$\wp_{1}'$ and $\wp_{2}'$ be as defined in (\ref{prtx10geneqns}) and (\ref{prtx11geneqns}) respectively.
The sets $D_{i}$, $i=1,4,8,10,11$ are chosen so that each one is the homomorphic image of a basis for $\textsf{W}_{1}$,$\textsf{W}_{4}$,$\textsf{W}_{8}$,$\textsf{W}_{10}\oplus\textsf{W}_{10}$,$\textsf{W}_{11}\oplus\textsf{W}_{11}$,  respectively,  under some $G$-isomorphism to $\pi_{{}_{\textsf{W}_{i}}}(I_{X_{1}}(2))$.
Then we have
\begin{align}
\nonumber\\
 D_{1}=&\{d_{1,7,1}\},\\  
&\mbox{ where } 
\qquad d_{1,7,1}=\pi_{{}_{\textsf{W}_{1}}}(y_{1}y_{2}).\nonumber\\
\nonumber\\
 D_{4}=&D_{4,1}\cup D_{4,2}, \\
&\mbox{ where }
\qquad D_{4,i}=\{d_{4,j,i}\mid   \mbox{$j$=1,$\dots$,6.} \},\mbox{ $i$=1,2,}\nonumber\\
& \mbox{   and } 
\qquad\:\: d_{4,j,i}= \rho(h_{7B})_{\zeta^j}(Q\cdot \vartheta_{i}).\nonumber\\
\nonumber\\
 D_{8}=&\{d_{8,j,1}\mid \mbox{$j$=1,$\dots$,7.} \}\cup\{d_{8,7,2} \},\\
&\mbox{ where }\nonumber\\
&\qquad\qquad d_{8,j,1}=\rho(h_{7B})_{\zeta^j}(Q\cdot q_{\zeta^2}) ,\nonumber\\
&\qquad\qquad d_{8,7,2}=\rho(h_{7B})_{1}(Q^2\cdot q_{\zeta^2}).\nonumber\\
\nonumber\\
D_{10}=&D_{10,1}\cup D_{10,2}, \\
&\mbox{ where }
\qquad D_{10,i}=\phi_{\theta_{1},\wp_{i}}(B_{1}),\nonumber\\
&\mbox{ $B_{1}$ is defined in formula (\ref{basB1}) and }\nonumber\\
&
\qquad \phi_{\theta_{1},\wp_{i}}:\mathbb{C}G_{<\theta_{1}>}\rightarrow\mathbb{C}G_{<\wp_{i}>}\nonumber
\end{align}
is the isomorphism of $G$-modules that maps $\theta_{1}$ to $\wp_{i}$. We write 
\begin{align} 
 d_{10,j,k}=\:&\phi_{\theta_{1},\wp_{1}}(\mbox{\small\ss{}}_{k,j}), \nonumber\\
d_{10,j,k+3}=&\phi_{\theta_{1},\wp_{2}}(\mbox{\small\ss{}}_{k,j}), k=1,2,3.\nonumber
\end{align}
One has for instance
\begin{align}
d_{10,2,1}=&\wp_{1},\nonumber\\
d_{10,2,k+1}=&\rho(h_{7B})_{\zeta^2}(Q\cdot d_{10,2,k}), \mbox{$k$=1,2.}\nonumber\\
d_{10,j,k}=&\rho(h_{7B})_{\zeta^j}(Q\cdot d_{10,2,k}), \mbox{$j$=1,3,4,5,6,7, $k$=1,2,3.}\nonumber
\end{align}
Similarly we define
\begin{align} 
 D_{11}=&D_{11,1}\cup D_{11,2}, \\
&\mbox{ where } \nonumber\\
&\qquad\qquad D_{11,i}=\{d_{11,j,3(i-1)+k}\mid \mbox{$j$=1,\dots,7, $k$=1,2,3.}  \},\mbox{$i$=1,2.}\nonumber\\
&\mbox{ and } \nonumber\\
&\qquad \qquad d_{11,j,3(i-1)+k}= \phi_{\theta_{1}',\wp_{i}'}(\mbox{\small\ss{}}_{k,j}').\nonumber
\end{align}

\nothing{
\begin{align} 
 &\quad\quad\{\pi_{{}_{\textsf{W}_{1}}}(y_{1}y_{2})\}\\
 &\quad\quad\{\rho(h_{7B})_{\zeta^j}(Q\cdot \vartheta_{i})\mid \mbox{j=1,...,6}\},\mbox{for $i=1,2$.} \\
&\quad\quad\{\rho(h_{7B})_{\zeta^j}(Q\cdot q_{\zeta^2})\mid \mbox{j=1,...,7}\}\cup\{\rho(h_{7B})_{1}(Q^2\cdot q_{\zeta^2})\}\\
D(x)&\;:=\{x,\,\rho(h_{7B})_{\zeta^2}\big(\mbox{\footnotesize$Q\cdot x$}\big),
\,\rho(h_{7B})_{\zeta^2}\Big( \mbox{\footnotesize$Q\cdot\rho(h_{7B})_{\zeta^2}\big(\mbox{\scriptsize$Q\cdot x$}\big)$}\Big)\},\\
&\quad\quad\mbox{for } x=\wp_{1},\wp_{2},\wp_{1}',\wp_{2}'.\nonumber\\
 E(x)&\;:=\{  \rho(h_{7B})_{\zeta^j}\big( Q\cdot d\big)\mid d \in D(x),j=1,3,4,5,6,7\}, \\
&\quad\quad\mbox{for } x=\wp_{1},\wp_{2},\wp_{1}',\wp_{2}'.\nonumber
\end{align}
}
\section{The space of cubic forms}\label{sectcubicforms}
We conclude the genus 17 case with  the following theorem.
\begin{theorem}\label{mainthm}
 The homogeneous ideal of the canonical embedding of $X_{1}$ in $\mathbb{P}^{16}$ is generated by the quadratic forms in it.
\end{theorem}

For the proof we used the basis defined in Section \ref{basI2} to show that  the space of cubic forms
generated by $I_{X_{1}}(2)$ is $I_{X_{1}}(3)$, namely  $\textsf{V}_{1}\cdot I_{X_{1}}(2)=I_{X_{1}}(3)$.
To do that it is enough to show that 
$\pi_{{}_{\textsf{W}_{i}}}(\textsf{V}_{1}\cdot I_{X_{1}}(2))=\pi_{{}_{\textsf{W}_{i}}}(I_{X_{1}}(3))$, $i=$1,\dots,11.
So we computed the dimensions for the eigenspaces 
 $\rho(h_{7B})_{\zeta^j}\Big(\pi_{{}_{\textsf{W}_{i}}}\big(\mbox{\scriptsize$\textsf{V}_{1}\cdot I_{X_{1}}(2)$}\big)\Big)$,
for some $j$. Thus the multiplicity of the representation $\textsf{W}_{i}$ in 
$\pi_{{}_{\textsf{W}_{i}}}(\textsf{V}_{1}\cdot I_{X_{1}}(2))$ can be deduced from the information in Table \ref{charpols}.
In Table \ref{jeigenspacesOfI3wi} we list the dimensions of the spaces 
$\rho(h_{7B})_{\zeta^j}\Big(\pi_{{}_{\textsf{W}_{i}}}\big(\mbox{\scriptsize$\textsf{V}_{1}\cdot I_{X_{1}}(2)$}\big)\Big)$
and  we also list bases for these spaces in order to help speed up verification since this calculation was very slow.
On the other hand, from Tables \ref{tabla:desh0KnX1} and \ref{tabla:desSnW9W2} we have
{\footnotesize
\begin{align}
H^0(X_{1},K_{X_{1}}^{\otimes 3})=&
\textsf{W}_{3}\oplus \textsf{W}_{5}\oplus\textsf{W}_{6}\oplus
\textsf{W}_{7}\oplus\textsf{W}_{9}\oplus\textsf{W}_{10}\oplus \textsf{W}_{11},\nonumber\\
\mbox{ and ~~~~~~~~~~~~~~}&\nonumber\\  
S^3(\textsf{V}_{1})\:=&
\textsf{W}_{2}^{4}\oplus \textsf{W}_{3}^{5}\oplus \textsf{W}_{4}\oplus  \textsf{W}_{5}^{8}\oplus 
\textsf{W}_{6}^{8}\oplus \textsf{W}_{7}^{8}\oplus \textsf{W}_{8}^{5}\oplus 
\textsf{W}_{9}^{13}\oplus \textsf{W}_{10}^{13}\oplus \textsf{W}_{11}^{13},\nonumber\\
\mbox{ then ~~~~~~~~~~~~~~}\nonumber\\  
H^0(\mathbb{P}^{16},\mathcal{I}_{X_{1}}(3))=&
\textsf{W}_{2}^{4}\oplus \textsf{W}_{3}^{4}\oplus \textsf{W}_{4}\oplus  \textsf{W}_{5}^{7}\oplus 
\textsf{W}_{6}^{7}\oplus \textsf{W}_{7}^{7}\oplus \textsf{W}_{8}^{5}\oplus 
\textsf{W}_{9}^{12}\oplus \textsf{W}_{10}^{12}\oplus \textsf{W}_{11}^{12}.
\end{align}
}
\begin{table}[!h]
\caption{Generators for $I_{X_{1}}(3)$\label{jeigenspacesOfI3wi}}
{\scriptsize
\begin{displaymath}
\begin{array}{l}
\begin{array}{lcc}
\hline\vspace{-.2cm}\\
\textsf{W}_{i},j&\begin{array}{c}\mbox{ dimension of }\\
                 \mbox{\scriptsize$\rho(h_{7B})_{\zeta^j}$}\!\Big(\pi_{{}_{\textsf{W}_{i}}}\!\big(\mbox{\scriptsize$\textsf{V}_{1}\!\!\cdot\!I_{X_{1}}(2)$}\big)\Big)\end{array}& \begin{array}{c}
 \mbox{ basis for } \\ \rho(h_{7B})_{\zeta^j}\Big(\pi_{{}_{\textsf{W}_{i}}}\big(\mbox{\scriptsize$\textsf{V}_{1}\cdot I_{X_{1}}(2)$}\big)\Big)\end{array}
                                                                                                                                              \\
\hline\vspace{-.2cm}\\
\textsf{W}_{1},-&0&--\\
\textsf{W}_{2},2&4&\pi_{{}_{\textsf{W}_{2}}}\!\left(\left\{\begin{array}{c}\ydgen{2}{3}{1}{7}{1},\ydgen{1}{3}{8}{1}{1},\\
                                                              \ydgen{0}{1}{10}{2}{1}, \ydgen{0}{1}{10}{2}{4} \end{array}\right\} \!\right)\\
\vspace{-7pt}\\
\textsf{W}_{3},3& 4& \pi_{{}_{\textsf{W}_{3}}}\!\left(\left\{\begin{array}{c} \ydgen{1}{3}{4}{2}{1},\ydgen{1}{3}{4}{2}{2},\\
                                  \ydgen{1}{1}{10}{2}{1}  ,\ydgen{1}{1}{10}{2}{4}    \end{array}\right\} \!\right)\\
\vspace{-7pt}\\
\textsf{W}_{4},2&1&\pi_{{}_{\textsf{W}_{4}}}\!\left(\left\{\ydgen{1}{3}{8}{1}{1}\right\} \!\right)\\
\vspace{-7pt}\\
\textsf{W}_{5},2& 7& \pi_{{}_{\textsf{W}_{5}}}\!\left(\left\{\begin{array}{c}\ydgen{1}{3}{4}{1}{1},\ydgen{1}{3}{4}{1}{2},\ydgen{1}{3}{8}{1}{1},\\
                                               \ydgen{1}{1}{10}{1}{1},  \ydgen{1}{1}{10}{1}{2},  \ydgen{1}{1}{10}{1}{4},\\
                                                    \ydgen{1}{1}{11}{1}{1}  \end{array}\right\} \!\right)\\
\vspace{-7pt}\\
\textsf{W}_{6},2&7&\pi_{{}_{\textsf{W}_{6}}}\!\left(\left\{\begin{array}{c} \ydgen{1}{1}{8}{1}{1},\ydgen{1}{1}{4}{1}{1},\ydgen{1}{1}{4}{1}{2},\\  
                                                                        \ydgen{1}{1}{10}{1}{1},   \ydgen{1}{1}{10}{1}{2}, \ydgen{1}{1}{10}{1}{4},\\  
                                                                   \ydgen{1}{1}{10}{1}{5}     \end{array}\right\} \!\right)\\                                     
\vspace{-7pt}\\
\textsf{W}_{7},2&7&\pi_{{}_{\textsf{W}_{7}}}\!\left(\left\{\begin{array}{c} \ydgen{1}{1}{8}{1}{1}, \ydgen{1}{1}{4}{1}{1}, \ydgen{1}{1}{4}{1}{2},\\
                                                                         \ydgen{1}{1}{10}{1}{1}, \ydgen{1}{1}{10}{1}{4}, \ydgen{1}{3}{10}{1}{1},\\
                                                                                      \ydgen{1}{3}{10}{1}{4} \end{array}\right\} \!\right)\\ 
\vspace{-7pt}\\
\textsf{W}_{8},2&5&\pi_{{}_{\textsf{W}_{8}}}\!\left(\left\{\begin{array}{c} \ydgen{1}{1}{10}{1}{1},\ydgen{1}{1}{10}{1}{2}, \ydgen{1}{1}{10}{1}{4},\\
                                                                             \ydgen{1}{1}{11}{1}{1}, \ydgen{1}{1}{11}{1}{2}  \end{array}\right\} \!\right)\\ 
\vspace{-7pt}\\
\textsf{W}_{9},2&24&\pi_{{}_{\textsf{W}_{9}}}\left(\begin{array}{c}
                                                                \ydDj{1}{1}{4}{1}\cup \ydDj{1}{2}{4}{1} \cup \ydDj{3}{1}{4}{6} \cup\\
                                                                     \ydDj{3}{2}{4}{6} \cup \ydDj{1}{1}{10}{1}\cup\\
\left\{\!\!\begin{array}{c}\ydgen{1}{1}{8}{1}{1},\ydgen{1}{2}{8}{1}{1},
\ydgen{1}{2}{10}{1}{1},\\ \ydgen{1}{2}{10}{1}{2}
\ydgen{1}{3}{10}{1}{1},\ydgen{1}{3}{10}{1}{2},\\
\ydgen{2}{1}{10}{7}{1},
\ydgen{1}{2}{10}{1}{4},\ydgen{1}{2}{10}{1}{5}\\
\ydgen{1}{3}{10}{1}{4}\end{array}\!\!\right\}
                      
                      \end{array}\right)\\
\vspace{-7pt}\\
\textsf{W}_{10},2&36&\pi_{{}_{\textsf{W}_{10}}}\!\left(\!\!\begin{array}{c}
\ydDj{1}{1}{4}{1}\cup\ydDj{1}{2}{4}{1}\cup\ydDj{3}{1}{4}{6}\cup\\
\ydDj{2}{1}{8}{7}\cup \ydDj{2}{2}{8}{7}\cup
\ydDj{1}{3}{11}{1}\cup\\
\ydDj{1}{1}{10,1}{1}\cup\ydDj{1}{2}{10,1}{1}\cup\\
\ydDj{2}{1}{10,1}{7}\cup \ydDj{2}{2}{10,1}{7}\cup\\

\left\{\!\!\begin{array}{c}\ydgen{1}{1}{8}{1}{1},\ydgen{1}{2}{8}{1}{1},
\ydgen{1}{3}{10}{1}{1},\\ \ydgen{1}{3}{10}{1}{2},
\ydgen{2}{3}{10}{7}{1},
\ydgen{1}{3}{10}{1}{4},\\ \ydgen{1}{3}{10}{1}{5},
\ydgen{2}{3}{10}{7}{4}
\end{array}\!\!\right\}
                                \end{array}\!\!\right)\\
\vspace{-7pt}\\
\textsf{W}_{11},2&36&\pi_{{}_{\textsf{W}_{11}}}\!\left(\!\!\!\begin{array}{c}
\ydDj{1}{1}{4}{1}\cup\ydDj{1}{2}{4}{1}\cup
\ydDj{3}{1}{4}{6}\cup  \\ \ydDj{2}{1}{8}{7}\cup
\ydDj{2}{2}{8}{7}\cup\ydDj{1}{3}{10}{1}\cup \\
\ydDj{1}{1}{10,1}{1}\cup\ydDj{1}{2}{10,1}{1}\cup\\
\ydDj{2}{1}{10,1}{7}\cup\ydDj{2}{2}{10,1}{7}\cup\\
\left\{\!\!\begin{array}{c}
\ydgen{1}{1}{8}{1}{1},\ydgen{1}{2}{8}{1}{1},
\ydgen{1}{3}{11}{1}{1}, \\ \ydgen{1}{3}{11}{1}{2},
\ydgen{2}{3}{11}{7}{2},
\ydgen{1}{3}{11}{1}{4},\\\ydgen{1}{3}{11}{1}{5},
\ydgen{2}{3}{11}{7}{4}
\end{array}\!\!\right\}
\end{array}\!\!\!\right)\\
\vspace{-3pt}\\\hline\\
\end{array}\\
\vspace{-2pt}
\mbox{Here we write $\left[y_{r}\right]_{\zeta^j}=\rho(h_{7B})_{\zeta^j}(y_{r})$ and }
\mbox{ $\left[ D_{*}\right]_{\zeta^{v}}=\left\{d_{i,j,k}\in D_{*}\mid j\equiv v \mod 7\right\}$}.
\end{array}
\end{displaymath}
}
\end{table}\\

\newpage
\section{The genus 14 case}\label{sectgenus14}
Let $X_{1}$ be the Hurwitz curve of genus 14 considered in \cite{hcg14}. Using the quadrics in the Theorem 2.5 of \cite{hcg14} one also proves the following.
\begin{theorem}\label{mainthm2}
 The homogeneous ideal of the canonical embedding of $X_{1}$ in $\mathbb{P}^{13}$ is generated quadrics.
\end{theorem}
\noindent We proceed exactly as in the proof of Theorem \ref{mainthm}. In this case we used the Table \ref{jeigenspacesOfI3wiGenus14} below, the Table 3 in \cite{hcg14} and the decomposition of $H^0(\mathbb{P}^{13},\mathcal{I}_{X_{1}}(3))$ into a sum of irreducible representations of $PSL(2,\mathbb{F}_{13})$. That is (use Tables 5 and 8 in \cite{hcg14}),
\begin{align}
H^0(\mathbb{P}^{13},\mathcal{I}_{X_{1}}(3))=&
\textsf{W}_{2}^{5}\oplus \textsf{W}_{3}^{5}\oplus \textsf{W}_{4}^{5}\oplus  \textsf{W}_{5}^{6}\oplus 
\textsf{W}_{6}^{5}\oplus \textsf{W}_{7}^{5}\oplus \textsf{W}_{8}^{4}\oplus 
\textsf{W}_{9}^{8}.
\end{align}
When computing the quadrics in Theorem 2.5 of \cite{hcg14} one should have into account two misprints: the first in the definition of $\epsilon_{i}$ (before formula (9)) where $\epsilon_{2}$ should be equal to $\rho(h_{7A})_{\zeta^{1}}(e_{3})$ and  the second  in the definition of $q$ (also before formula (9)) where $q=\rho(h_{7A})_{\zeta^{5}}(\pi_{{}_{\textsf{W}_{6}}(y_{1}^{2})})$ should be $q=\rho(h_{7A})_{\zeta^{5}}(\pi_{{}_{\textsf{W}_{6}}(y_{2}^{2})})$.
\begin{table}[!h]
\caption{Generators for $I_{X_{1}}(3)$ in the genus 14 case.\label{jeigenspacesOfI3wiGenus14}}
{\scriptsize
\begin{displaymath}
\begin{array}{l}
  \begin{array}{lcc}
   \hline\vspace{-.2cm}\\
    \textsf{W}_{i},j & {\begin{array}{c}\mbox{ dimension of }\\
                 \mbox{\scriptsize$\rho(h_{13A})_{\nu^j}$}\!\Big(\pi_{{}_{\textsf{W}_{i}}}\!\big(\mbox{\scriptsize$\textsf{V}_{1}\!\!\cdot\! I_{X_{1}}(2)$}\big)\Big)
                        \end{array}} &           \begin{array}{c}
 \mbox{ basis for } \\ \rho(h_{13A})_{\nu^j}\Big(\pi_{{}_{\textsf{W}_{i}}}\big(\mbox{\scriptsize$\textsf{V}_{1}\cdot I_{X_{1}}(2)$}\big)\Big)\end{array}                                \\
  \hline\vspace{-.2cm}\\

\textsf{W}_{1},-&0&--\\
\textsf{W}_{2},2&5&\pi_{{}_{\textsf{W}_{2}}}\!\left(\left\{\begin{array}{c}\ydgenfrstp{3}{1}{5}{12}{1},\ydgenfrstp{1}{1}{6}{1}{1},\\
                                                              \ydgenfrstp{1}{1}{7}{1}{1}, \ydgenfrstp{1}{1}{8}{1}{1},\ydgenfrstp{1}{1}{8}{1}{3} \end{array}\right\} \!\right)\\
\vspace{-7pt}\\
\textsf{W}_{3},1& 5& \pi_{{}_{\textsf{W}_{3}}}\!\left(\left\{\begin{array}{c} \ydgenfrstp{2}{1}{5}{12}{1},\ydgenfrstp{2}{1}{6}{12}{1},\\
                                  \ydgenfrstp{1}{1}{7}{13}{1}  ,\ydgenfrstp{1}{1}{8}{13}{1}   \ydgenfrstp{1}{1}{8}{13}{2}  \end{array}\right\} \!\right)\\
\vspace{-7pt}\\
\textsf{W}_{4},1&5&\pi_{{}_{\textsf{W}_{4}}}\!\left(\left\{\begin{array}{c}\ydgenfrstp{2}{1}{5}{12}{1},
\ydgenfrstp{3}{1}{5}{11}{1},\\
\ydgenfrstp{2}{1}{6}{12}{1}, \ydgenfrstp{3}{1}{6}{11}{1}\ydgenfrstp{1}{1}{7}{13}{1} \end{array}\right\} \!\right)\\
\vspace{-7pt}\\
\textsf{W}_{5},1& 6& \pi_{{}_{\textsf{W}_{5}}}\!\left(\left\{\begin{array}{c}\ydgenfrstp{3}{1}{5}{11}{1},\ydgenfrstp{4}{1}{5}{10}{1},\ydgenfrstp{2}{1}{6}{12}{1},\\
                                               \ydgenfrstp{3}{1}{6}{11}{1},  \ydgenfrstp{1}{1}{7}{13}{1},  \ydgenfrstp{2}{1}{7}{12}{1}
                                           \end{array}\right\} \!\right)\\
\vspace{-7pt}\\
\textsf{W}_{6},1&5&\pi_{{}_{\textsf{W}_{6}}}\!\left(\left\{\begin{array}{c} \ydgenfrstp{2}{1}{5}{12}{1},\ydgenfrstp{3}{1}{5}{11}{1},\ydgenfrstp{3}{1}{6}{11}{1},\\  
                                                                        \ydgenfrstp{4}{1}{6}{10}{1},   \ydgenfrstp{1}{1}{7}{13}{1} 
                                                                  \end{array}\right\} \!\right)\\                                     
\vspace{-7pt}\\
\textsf{W}_{7},1&5&\pi_{{}_{\textsf{W}_{7}}}\!\left(\left\{\begin{array}{c} \ydgenfrstp{2}{1}{5}{12}{1}, \ydgenfrstp{3}{1}{5}{11}{1}, \ydgenfrstp{2}{1}{6}{12}{1},\\
                                                                         \ydgenfrstp{3}{1}{6}{11}{1}, \ydgenfrstp{1}{1}{7}{13}{1}\end{array}\right\} \!\right)\\ 
\vspace{-7pt}\\
\textsf{W}_{8},1&4&\pi_{{}_{\textsf{W}_{8}}}\!\left(\left\{\begin{array}{c} \ydgenfrstp{2}{1}{5}{12}{1}, \ydgenfrstp{3}{1}{5}{11}{1}, \ydgenfrstp{2}{1}{6}{12}{1},\\
                                                                         \ydgenfrstp{3}{1}{6}{11}{1} \end{array}\right\} \!\right)\\ 
\vspace{-7pt}\\
\textsf{W}_{9},1&8&\pi_{{}_{\textsf{W}_{9}}}\!\left(\left\{\begin{array}{c} \ydgenfrstp{1}{1}{8}{13}{1}, \ydgenfrstp{1}{1}{8}{13}{2}, \ydgenfrstp{2}{1}{8}{12}{1},\\
                                                                         \ydgenfrstp{3}{1}{7}{11}{1}, \ydgenfrstp{2}{1}{5}{12}{1}, \ydgenfrstp{3}{1}{5}{11}{1},\\ 
\ydgenfrstp{2}{1}{6}{12}{1}, \ydgenfrstp{3}{1}{6}{11}{1} \end{array}\right\} \!\right)\\ 
\vspace{-3pt}\\\hline\\
  \end{array}\\
\vspace{-2pt}
\mbox{Using the notation of \cite{hcg14} we write $\textsf{V}_{1}=<y_{1},\dots,y_{14}>$},\\
\mbox{$\left[y_{r}\right]_{\nu^j}=\rho(h_{13A})_{\nu^j}(y_{r})$, $\nu=e^{2i\pi/13}$},\\
\mbox{$d_{i,j,1}= \rho(h_{13A})_{\nu^{j}}(\pi_{\textsf{W}_{5}}(y_{2}^{2}))$ for $i=5,6$ and $j=1\dots12$},\\
\mbox{$d_{7,j,1}= \rho(h_{13A})_{\nu^{j}}(h_{7A}\cdot q_{0})$ for  $j=1\dots13$},\\
\mbox{$d_{8,j,2k-1}=\rho(h_{13A})_{\nu^{j}}(h_{7A}\cdot \vartheta_{k})$ for $k=1,2$ and  $j=1\dots 13$},\\
\mbox{$d_{8,j,2k}=\rho(h_{13A})_{\nu^{j}}(h_{7A}^{2}\cdot \vartheta_{k})$ for $k=1,2$ and $j=13$}.\\
\end{array}
\end{displaymath}
}
\end{table}

~\\

\newpage
\noindent{\bf Acknowledgments.} Some calculations were carried out at the KanBalam  supercomputer of DGTIC--UNAM.\\

{\setlength{\baselineskip}%
           {0.1\baselineskip}
\noindent Israel Moreno Mej\'{\i}a.\\
Instituto de Matem\'aticas \\
Universidad Nacional Aut\'onoma de M\'exico\\
M\'exico, D.F. C.P. 04510\\
M\'exico.
}

\vspace{.5cm}
\noindent E-mail address: {\bf israel@math.unam.mx}

\end{document}